\documentclass[12pt]{amsart}
\usepackage{amssymb,latexsym}
\newtheorem{theorem}{Theorem}
\newtheorem{lemma}[theorem]{Lemma}

\def\ZZ{{\mathbb Z}}
\def\RR{{\mathbb R}}

\def\se{{\subseteq}}

\def\QQ{{\mathbb Q}}

\begin{document}

\title[isoperimetric  inequalities for arithmetic groups]{Exponential higher dimensional isoperimetric  inequalities for some arithmetic groups}

\author{Kevin Wortman}

\begin{abstract}

We show that arithmetic subgroups of semisimple groups of relative $\QQ$-type $A_n$, $B_n$, $C_n$, $D_n$, $E_6$, or $E_7$ have an exponential lower bound to their isoperimetric inequality in the dimension that is $1$ less than the real rank of the semisimple group.

\end{abstract}

\maketitle

Let $\bf{G}$ be a connected, semisimple, $\QQ$-group that is almost simple over $\QQ$. Let $X$ be the symmetric space of noncompact type associated with ${\bf G}(\RR)$ and let  $X_\ZZ $ be a contractible subspace of $X$ that is a finite Hausdorff distance from some ${\bf G}(\ZZ)$-orbit in $X$; Raghunathan proved that such a space exists \cite{Ra 1}. We denote the $\RR$-rank of ${\bf G}$ by ${\text{rk}}_\RR{\bf G}$.

Given a homology $n$-cycle $Y\se X_\ZZ$ we let $v_X(Y)$ be the infimum of the volumes of all $(n+1)$-chains $ B \se X$ such that $\partial B=Y$. Similarly, we let $v_\ZZ(Y)$ be the infimum of the volumes of all $(n+1)$-chains $ B \se X_\ZZ$ such that $\partial B=Y$. We define the ratio $$ R_n(Y) = \frac{v_\ZZ(Y)}{v_X(Y)}$$ and we let $R_n({\bf G}(\ZZ)) : \RR _{>0} \rightarrow \RR_{\geq 1}$ be the function $$R_n({\bf G}(\ZZ))(L)=\sup \{\, R_n(Y) \mid {\text{vol}}(Y) \leq L\,\}$$
These functions measure a contrast between the geometries of ${\bf{G}}(\ZZ)$ and $X$.

 Clearly if $\bf G$ is $\QQ$-anisotropic (or equivalently, if  ${\bf G}(\ZZ)$ is cocompact in ${\bf G}(\RR)$) then we may take  $X_\ZZ = X$ so that $R_n({\bf G}(\ZZ) )=1$ for all $n$. 

The case is different when $\bf G$ is $\QQ$-isotropic, or equivalently, if  ${\bf G}(\ZZ)$ is non-cocompact in ${\bf G}(\RR)$.

Leuzinger-Pittet conjectured that $R_{{\text{{rk}}}_\RR{\bf G}-1}({\bf G}(\ZZ))$ is bounded below by an exponential when $\bf G$ is $\QQ$-isotropic  \cite{L-P}. The conjecture in the case ${\text{rk}}_\RR{\bf G}=1$ is equivalent to the well-known observation that the word metric for non-cocompact lattices in rank one real simple Lie groups is exponentially distorted in its  corresponding symmetric space. Prior to \cite{L-P}, the conjecture was evidenced by other authors in some cases. It was proved by Epstein-Thurston when ${\bf G}(\ZZ)={\bf SL_k}(\ZZ)$  \cite{Ep}, by Pittet when ${\bf G}(\ZZ)={\bf SL_2}(\mathcal{O})$ and $\mathcal{O}$ is a ring of integers in a totally real number field  \cite{Pi}, by Hattori when ${\bf G}(\ZZ)={\bf SL_k}(\mathcal{O})$ and $\mathcal{O}$ is a ring of integers in a totally real number field \cite{Ha 1}, and by Leuzinger-Pittet when ${\text{rk}}_\RR{\bf G}=2$  \cite{L-P}.

This paper contributes to the verification of the Leuzinger-Pittet conjecture by proving

\begin{theorem}\label{t:m}
Let $\bf{G}$ be as in the introductory paragraph and assume that $\bf G$ is $\QQ$-isotropic. Furthermore, suppose the $\QQ$-relative root system of $\bf G$ is of type $A_n$, $B_n$, $C_n$, $D_n$, $E_6$, or $E_7$. Then there exist constants $C>0$ and $L_0>0$ such that $$R_{{\text{ \emph{rk}}}_\RR{\bf G}-1}({\bf G}(\ZZ))(L) \geq e^{CL}$$ for any $L> L_0 $.
\end{theorem}

\subsection{Example.} Let $\mathcal{O}$ be the ring of integers in a number field $K$, and let ${\bf G}={\bf R}_{K/\QQ}{\bf SL_k}$ where ${\bf R}_{K/\QQ}$ is the restriction of scalars functor. Then ${\bf G}(\ZZ )={\bf SL_k}(\mathcal{O})$, $\bf G$ is $\QQ$-isotropic, $\bf G$ has a $\QQ$-relative root system of type $A_{k-1}$, and ${\text{rk}}_\RR{\bf G}=(k-1)S$ where $S$ is the number of inequivalent archimedean valuations on $K$. Therefore, $R_{(k-1)S-1}({\bf SL_k}(\mathcal{O}))$ is bounded below by an exponential.

\subsection{Non-nonpositive curvature of arithemtic groups.} If ${\bf G}(\ZZ)$ satisfied a reasonable notion of nonpositive curvature (including $\text{CAT}(0)$ or combable, for example), we would expect polynomial bounds on isoperimetric inequalities for ${\bf G}(\ZZ)$. Thus, not only does Theorem~\ref{t:m} provide a measure of the difference between ${\bf G}(\ZZ)$ and $X$, it also exhibits \emph{non}-nonpositive curvature tendencies for ${\bf G}(\ZZ)$ when $\bf G$ is $\QQ$-isotropic and ${\text{rk}}_\RR{\bf G}>1$.

\subsection{Type restriction.} Our proof of Theorem~\ref{t:m} excludes the remaining types -- $G_2$, $F_4$, $E_8$, and $BC_n$ -- because groups of these types do not contain proper parabolic subgroups whose unipotent radicals are abelian. Our techniques require an abelian unipotent radical of a maximal $\QQ$-parabolic subgroup of $\bf G$ to construct cycles in $X_\ZZ$.

\subsection{Related results.} It is an open question whether $R_{n}({\bf G}(\ZZ))$ is bounded above by a constant when $n<{\text{rk}}_\RR{\bf G}-1$. When $n=0$ it is; this is a theorem of Lubotzky-Mozes-Raghunathan \cite{L-M-R}.

Dru\c{t}u showed that if the $\QQ$-relative root system of $\bf G$ is of type $A_1$ or $BC_1$, then for any $\varepsilon >0$, ${\bf G}(\ZZ)$ has a Dehn function that is bounded above by $L^{2+\varepsilon}$ for $L$ sufficiently large \cite{Dr}.

Young proved that ${\bf{G}}(\ZZ)={\bf{SL_k}(\ZZ)}$ has a quadratic Dehn function if $k\geq 5$ \cite{Yo}.

Gromov proved that all of the functions $R_{n}({\bf G}(\ZZ))$ are bounded above by an exponential function, and Leuzinger later provided a more detailed proof of this fact ($5.A_7$ \cite{Gr} and Corollary 5.4 \cite{Le}). 

\section{Choice of parabolic}

Let $\bf T \leq G$ be a maximal $\QQ$-split torus in $\bf G$. We let $\Phi _\QQ$ be the roots of $\bf G$ with respect to $\bf T$. Choose an ordering on $\Phi _\QQ$. We denote the corresponding sets of simple and positive roots by  $\Delta _\QQ$ and $\Phi _\QQ^+$ respectively.

If $I \se \Delta _\QQ$, we let $[I] \se \Phi _\QQ$ be the set of roots that are linear combinations of elements in $I$, and we let $\Phi _\QQ (I)^+ =\Phi _\QQ^+ - [I]$.

For each $\alpha \in  \Phi _\QQ$, we let ${\bf U}_{\alpha} \leq {\bf G}$ be the root subgroup associated with $\alpha$. For $J \se\Phi _\QQ$, we let ${\bf U}_J=\prod _{\alpha \in J}{\bf U_{\alpha}}$.

We define ${\bf T}_I = \cap _{\alpha \in I}{\text{Ker}}(\alpha )^\circ$  where the superscript $\circ$ denotes the connected component of the identity, and we label the centralizer of ${\bf T}_I$ in $\bf G$ by ${\bf Z_G(T}_I)$.

\subsection{Maximal parabolics with abelian unipotent radicals.} For any $\alpha _0 \in \Delta _\QQ$, we let ${\bf P}_{\alpha_0}$ be the maximal proper parabolic subgroup of $\bf{G}$ given by $ {\bf U}_{\Phi _\QQ (\Delta _\QQ - \alpha _0)^+}{\bf Z_G(T}_{\Delta _\QQ - \alpha _0})   $. The unipotent radical of ${\bf P}_{\alpha_0}$ is  ${\bf U}_{\Phi _\QQ (\Delta _\QQ - \alpha _0)^+}$.

\begin{lemma}\label{l:fa} There is some $\alpha _0 \in \Delta _\QQ$ such that ${\bf U}_{\Phi _\QQ (\Delta _\QQ - \alpha _0)^+}$ is abelian.
\end{lemma}

\begin{proof}
Suppose $\Delta _\QQ=\{\alpha_1, \alpha_2,...,\alpha_{k}\}$. The set of positive roots $\Phi ^+ _\QQ$ contains a  ``highest root" $\sum_i n_i\alpha _i$ for positive integers $n_i$ such that if  $\sum_i m_i\alpha _i \in \Phi ^+ _\QQ$, then $m_i\leq n_i$ (\cite{Bo}, VI 1 8).

Given that $\Phi _\QQ$ is a root system of type  $A_n$, $B_n$, $C_n$, $D_n$, $E_6$, or $E_7$, there is some $\alpha _0 \in \{\alpha_1, \alpha_2,...,\alpha_{k}\}$ such that $n_0=1$; consult the list of root systems in the appendix of \cite{Bo}.

Since any $\sum_i m_i\alpha _i \in \Phi _\QQ (\Delta _\QQ - \alpha _0)^+$ has $m_0>0$, it follows that any $\sum_i m_i\alpha _i \in \Phi _\QQ (\Delta _\QQ - \alpha _0)^+$ has $m_0=1$, and thus the sum of two elements in $ \Phi _\QQ (\Delta _\QQ - \alpha _0)^+$ is not a root.

Therefore, given $\tau _1,\tau_2 \in  \Phi _\QQ (\Delta _\QQ - \alpha _0)^+$, we have $$[{\bf U}_{\tau_1},{\bf U}_{\tau _2}] \se {\bf U}_{\tau _1 + \tau _2}=1$$

\end{proof}

In what remains, we let ${\bf P=P}_{\alpha _0}$, we let $U_P$ be the real points of ${\bf U}_{\Phi _\QQ (\Delta _\QQ - \alpha _0)^+}$. Thus, we can rephrase Lemma~\ref{l:fa} as 

\begin{lemma}\label{l:sa} $U_P$ is abelian.
\end{lemma}

\subsection{A contracting ray.} Recall that ${\bf T}_{\Delta _\QQ - \alpha _0} \leq {\bf Z_G(T}_{\Delta _\QQ - \alpha _0}) \leq {\bf P} $ is a $1$-dimensional $\QQ$-split torus.
Choose $a_+ \in {\bf{T}}_ {\Delta _\QQ - \alpha _0}(\RR) $ such that $\alpha _0 (a_+)>1$ and such that the distance in  ${\bf{T}}_ {\Delta _\QQ - \alpha _0}(\RR)$ between $1$ and $a_+$ equals $1$.

The Lie algebra of $U_P$ is $\mathfrak{u}$.

\begin{lemma}\label{l:ad}
There is some $s>0$ such that for any  $v \in \mathfrak{u}$ $$Ad(a_+^{t})v  = e^{st} v $$
\end{lemma}

\begin{proof}

Recall that  $$\mathfrak{u}=\prod _{\beta \in \Phi _\QQ (\Delta _\QQ - \alpha _0)^+ }\mathfrak{u}_\beta$$ where $$\mathfrak{u_\beta}=\{\, v \in \mathfrak{u} \mid Ad(x)v=\beta(x)v{\text{ for all }} x\in {\bf T} \,\}$$ 

If $\beta \in \Phi _\QQ (\Delta _\QQ - \alpha _0)^+$, then $\beta = \alpha _0 +\sum _{\alpha _i \in \Delta _\QQ - \alpha _0} n_i \alpha _i$. Since $a_+ \in \cap _{\alpha_i \in \Delta _\QQ - \alpha _0 }{\text{Ker}}(\alpha_i )^\circ$, we have $\beta (a_+)=\alpha _0(a_+)$ and thus for $v \in \mathfrak{u}$, it follows that $Ad(a_+)v=\alpha_0(a_+)v$. 

Let $s=\log \big( \alpha _0 (a_+) \big)$.

\end{proof}

\section{A horoball in the symmetric space, disjoint from $X_\ZZ$}

\begin{lemma} There is a maximal $\QQ$-torus $\bf A \leq G$ such that the maximal $\QQ$-split torus of $\bf A$ is ${\bf T}_{\Delta _\QQ - \alpha _0}$ and such that $\bf A$ contains a maximal $\RR$-split torus of $\bf G$.
\end{lemma}

\begin{proof}
See Proposition 3.3 in \cite{B-W} where $K=\QQ$, $\mathbf{H}=\mathbf{G}$, $\mathbf{T_1}={\bf T}_{\Delta _\QQ - \alpha _0}$, $S=\{v\}$, and $K_v=\RR$. \end{proof}

Let $\bf Q$ be a minimal parabolic that contains $\bf A$ and is contained in $\bf P$. We let $\Phi _\RR$ be the roots of $\bf G$ with respect to the maximal $\RR$-split subtorus of $\mathbf{A}$,  $\Delta _\RR$ be the collection of simple roots given by $\bf Q$, and $\Phi_\RR^+$ be the corresponding positive roots.

\subsection{Alternate descriptions of the symmetric space.} Let $G={\bf G}(\RR)$ and let $A \leq G$ be the $\RR$-points of the maximal $\RR$-split subtorus of ${\bf A}$. Recall that ${\bf A}(\RR)=AB$ for some compact group $B \leq {\bf A}(\RR)$.

Choose a maximal compact subgroup $K \leq G$ that contains $B$. Then $G/K$ is a symmetric space that $G$ acts on by isometries. We name this symmetric space $X$.

Let $U_Q$ be the group of real points of the unipotent radical of $\bf Q$. By the Iwasawa decomposition, $U_Q A$ acts simply transitively on $X$ and we identify $X$ with  $U_Q A$. In this description of $X$, $A$ is a flat.

 \subsection{Integral translations in a flat.} By Dirichlet's units theorem, ${\bf A}(\ZZ) $ contains a finite index free abelian subgroup of rank ${\text{rk}}_\RR({\bf G})-1={\text{dim}}(A)-1$. Thus, if $A_\ZZ$ is the convex hull in $X$ of the ${\bf A}(\ZZ) $-orbit of the point $1 \in U_QA =X$, then $A_\ZZ$ is a codimension-1 Euclidean subspace of the flat $A$, and ${\bf A}(\ZZ) $ acts cocompactly on $A_\ZZ$.
 We may assume $A_\ZZ \se X_\ZZ$.

\subsection{Horoballs.} Notice that $\{a_+^t\}_{t>0}$ defines a unit-speed geodesic ray that limits to a point in $A ^\infty$ which we denote $a_+^\infty$.
We let $b_{a_+^t} : U_QA \rightarrow \RR$ be the Busemann function corresponding to the geodesic ray   $\{a_+^t\}_{t>0}$. That is, $$b_{a_+^t}(x)=\lim _{t\to \infty}[d(x,a_+^t)-t]$$

We let $A_0 \leq A$ be the codimension-$1$ subspace of $A$ consisting of those $a \in A$ for which $b_{a_+^t}(a)=0$. Thus, $A_0$ is orthogonal to $a_+^\RR$.

\begin{lemma}\label{l:b}
For $T\in \RR$, $(b_{a_+^t})^{-1}(-T)=U_QA_0a_+^T$.
\end{lemma}

\begin{proof}
We first show that for $u \in U_Q$ and $x \in X$, $b_{a_+^t}(x)=b_{ua_+^t}(x)$ Where $b_{ua_+^t}$ is the Busemann function for the ray  $\{ua_+^t\}_{t>0}$ .

Notice that $U_Q=U_PU_a$ where  $U_a \leq {\bf Z_G(T}_{\Delta _\QQ - \alpha _0}) (\RR)$ is a unipotent group whose elements commute with $a_+$.

If $u \in U_P$, then Lemma~\ref{l:ad} implies 
 $$d(a_+^{t}, ua_+^t)=d(1,a_+^{-t}ua_+^t) \to 0$$ Therefore,
$$ b_{a_+^t}(x)=\lim _{t\to \infty}[d(x,a_+^t)-t]=\lim _{t\to \infty}[d(x,ua_+^t)-t]=b_{ua_+^t}(x)  $$

The quotient map of a Lie group by a normal subgroup is distance nonincreasing. Because $U_P$ is normal in $U_Q A$, and because $a_+^\RR$ is normal in $U_aA$, the following composition is distance nonincreasing $$U_Q A \rightarrow U_aA \rightarrow U_aA_0  $$ 

We denote the geodesic between points $z,w \in X$ by $\overline{z,w}$. Orthogonality of $A_0$ and $a_+^\RR$ and the conclusion of the above paragraph show that for any $u \in U_a$, $\overline{1,u}$ is orthogonal to $a_+^\RR$ at $1$ and to $ua_+^\RR$ at $u$ and thus that $\overline{a_+^t,ua_+^t}$ is orthogonal to $a_+^\RR$ at $a_+^t$ and to $ua_+^\RR$ at $ua_+^t$. Furthermore, the length of $\overline{a_+^t,ua_+^t}$ is independent of $t$ since $u$ commutes with $a_+^t$.

Notice that the angle between $\overline{a_+^t,x}$ and $\overline{a_+^t,1}$ limits to $0$ as $t\to \infty$. Similarly, the angle between $\overline{ua_+^t, x}$ and $\overline{ua_+^t, u}$ limits to $0$. Hence, the triangle in $X$ with vertices $a_+^t$, $ua_+^t$, and $x$ approaches a triangle with angles $\frac\pi2$, $\frac\pi2$, and $0$. That is $$d(x,a_+^t)-d(x,ua_+^t) \to 0$$
Consequently, for $u \in U_a$ we have 
$$ b_{a_+^t}(x)=\lim _{t\to \infty}[d(x,a_+^t)-t]=\lim _{t\to \infty}[d(x,ua_+^t)-t]=b_{ua_+^t}(x)  $$

Therefore, for $u \in U_Q$,  $b_{a_+^t}(u^{-1} x)=b_{ua_+^t}(x)=b_{a_+^t}(x)$, and it follows that  $U_Q(b_{a_+^t})^{-1}(-T)=(b_{a_+^t})^{-1}(-T)$.
The lemma is a combination of this last fact together with $A_0 a_+^T \se (b_{a_+^t})^{-1}(-T)$.

\end{proof}

\begin{lemma}\label{l:T}
For some $T>0$, $X_\ZZ \se U_QA_0a_+^{(-\infty,T]}$.

\end{lemma}

\begin{proof}

Theorem A of \cite{Ha 2} states that $X_\ZZ \se (b_{a_+^t})^{-1}[-T,\infty) $ for some $T>0$, and $(b_{a_+^t})^{-1}[-T,\infty)=U_QA_0a_+^{(-\infty,T]}$ by Lemma~\ref{l:b}.

\end{proof}

\subsection{Projecting onto a horosphere.} Let $\pi: U_QA_0a_+^\RR \rightarrow U_QA_0$ be the obvious projection of $X$ onto the horosphere  $(b_{a_+^t})^{-1}(0)$.

\begin{lemma}\label{l:push} There is some $M >0$ such that for any $x_1,x_2 \in X_\ZZ$, we have $d(x_1,x_2) +M \geq d(\pi (x_1),\pi(x_2))$.

\end{lemma}

\begin{proof}
Recall that $U_Q=U_PU_a$ where elements of $U_a \leq \bf P $, and elements of $A_0$, commute with $a_+$. Similar to Lemma~\ref{l:ad}, we have that for any $t>0$ and any $v$ in the Lie algebra of $U_QA$ that 
$$||Ad(a_+^{-t})v || \leq ||v|| $$

Let $T$ be as in Lemma~\ref{l:T} and define $\pi _T :U_QA_0a_+^{(-\infty, T]} \rightarrow U_QA_0a_+^T$ by
 $\pi _T = R_{a_+^T} \circ \pi$ where  $R_{a_+^T}$ is right multiplication by  $a_+^T$.

We claim that $\pi _T $ is distance nonincreasing. To see this, first let $v$ be a tangent vector to $X$ at the point $a_+^t$ for some $t \leq T$. With $|| \cdot ||_x$ as the norm at $x$, and $f_*$ as the differential of $f$, we have

\begin{align*}
|| (\pi _T)_*v||_{\pi _T(a_+^t)} &  = || (R_{a_+^{T-t}})_*v||_{a_+^T} \\
& = || (L_{a_+^{t-T}})_*(R_{a_+^{T-t}})_*v||_{a_+^t} \\
& = || Ad(a_+^{t-T}) v||_{a_+^t} \\
& \leq || v ||_{a_+^t} 
\end{align*}

Left-translations by $U_QA_0$ show that for any $x \in U_QA_0a_+^{(-\infty, T]} $, and any $v \in T_xX$, $$ || (\pi _T)_*v||_{\pi_T(x)} \leq ||v||_{x}  $$

For any path $c :[0,1] \rightarrow U_QA_0a_+^{(-\infty, T]} $, apply $\pi _T$ to those segments contained in $U_QA_0a_+^{(-\infty, T]} $ to define a path between $\pi_T(c(0))$ and $\pi_T(c(1))$. This new path will have its length bounded above by the length of $c$ as is easily verified from the inequality on norms of vectors from above. This confirms our claim that $\pi_T$ is distance nonincreasing.

To confirm the lemma, notice that similarly, the map $R_{a_+^{-T}} : U_QA_0a_+^T \rightarrow U_QA_0$ translates all point in $X$ a distance of
$$d(x,  R_{a_+^{-T}}(x))=d(1,a_+^{-T})$$
Therefore, $$d(R_{a_+^{-T}}(x_1), R_{a_+^{-T}}(x_2)) \leq d(x_1,x_2)+2d(1,a_+^{-T})$$
The lemma follows as $\pi = R_{a_+^{-T}} \circ \pi_T$.
\end{proof}

\section{Choice of a cell in $X_\ZZ$}

We want to construct a cycle $Y \se X_\ZZ$. In this section we begin by constructing a cell $F \se A_0$ that will be used in the construction of $Y$.

\begin{lemma}\label{l:now}
$A_0 \se X_\ZZ$.
\end{lemma}

\begin{proof}

Both $A_0$ and the convex hull of $A_\ZZ$ are codimension $1$ subspaces of $A$. Since $A_\ZZ \se X_\ZZ \se U_QA_0a_+^{(-\infty,T)}$ we have that $A_\ZZ \se  A_0a_+^{(-\infty,T)}$. Therefore $A_\ZZ$ and $A_0$ are parallel hyperplanes. Since the both contain $1$, they are equal.

\end{proof}

 Let $X^\infty$ be the spherical Tits building  for $X=U_QA$, and let $A ^\infty \se X^\infty$ be the apartment given by $A$. Let $\Pi ^\infty  \se X^\infty$ be the simplex given by ${\bf P}$ and let $\Pi_- ^\infty  \se X^\infty$ be the simplex opposite of $\Pi ^\infty $ in $A ^\infty$, or equivalently, $\Pi_- ^\infty$ is the simplex given by the parabolic group ${\bf P^-}= {\bf U}_{\Phi _\QQ (\Delta _\QQ - \alpha _0)^-}{\bf Z_G(T}_{\Delta _\QQ - \alpha _0})   $.

Denote  the star of  $\Pi _-^\infty $ in $A ^\infty$ by $\Sigma \se A ^\infty$. Note that $\Sigma$ is homeomorphic to a ${\text{rk}}_\RR({\bf G})-1$ ball. We denote the codimension $1$ faces of $\Sigma$ as $\Sigma_1,...,\Sigma_n$.

\subsection{$A_0^\infty$ and $\Sigma$ are disjoint.} Let $\Psi \subseteq \Phi _\RR$ be such that ${\bf U}_\Psi={\bf R_u(P^-)}$. Given $b \in A_0$ we define the following sets of roots:
$$C(b)=\{\,\beta \in  \Psi \mid \beta (b) >1\,\}$$ 
$$Z(b)=\{\,\beta \in  \Psi \mid \beta (b) =1\,\}$$
$$E(b)=\{\,\beta \in  \Psi \mid \beta (b) <1\,\}$$

Thus, if $U_{C(b)}$ are the real points of ${\bf U}_{C(b)}$ etc., then ${\bf R_u(P^-)}(\RR)=U_{C(b)}U_{Z(b)}U_{E(b)}$.

\begin{lemma}\label{l:sweep} There is a sequence $\gamma _n \in {\bf R_u(P^-)}(\ZZ)-1$ such that $d(\gamma _n, U_{C(b)}) \to 0$.
\end{lemma}

\begin{proof} There is a $\QQ$-isomorphism of the variety ${\bf R_u(P^-)}$ with affine space that maps ${\bf U}_{C(b)}$ onto an affine subspace. Therefore, the problem reduces to showing that the distance between $\ZZ^n -1$ and a line in $\RR^n$ that passes through the origin is bounded above by any positive number, and this is well known.

\end{proof}

\begin{lemma}\label{l:disinf} $A_0^\infty \cap \Sigma = \emptyset$

\end{lemma}

\begin{proof}
Suppose $A_0^\infty \cap \Sigma \neq \emptyset$. Then there is some $b \in A_0$ such that $b^\infty \in \Sigma $ where $b^\infty = \lim _{t \to \infty}b^t$.

If $\mathfrak{C} \se \Sigma$ is a chamber, then $\Pi _-^\infty \se \mathfrak{C}$. Hence, the minimal $\RR$-parabolic subgroup corresponding to $\mathfrak{C}$ contains $\bf{R_u(P^-)}$ and thus elements of ${\bf R_u(P^-)}(\RR)$ fix $\mathfrak{C}$ pointwise. That is, elements of ${\bf R_u(P^-)}(\RR)$ fix $\Sigma$ pointwise, so they fix $b^\infty$.

Let $u \in {\bf R_u(P^-)}(\RR)$. Then $u b^\infty =b^\infty$, so $d(u b^t,b^t)$ is bounded, so $\{b^{-t}u b^t\}_{t>0}$ is bounded. It follows that $\beta (b^{-1})\leq 1$ for all $\beta \in \Psi$, or equivalently that $\beta (b)\geq 1$. Hence, $E(b) =\emptyset$ and ${\bf R_u(P^-)}(\RR)=U_{C(b)}U_{Z(b)}$.

Now we use Lemma~\ref{l:sweep}. For any $n \in \mathbb{N}$, there exists $\gamma_n \in {\bf R_u(P^-)}(\ZZ)-1$ with $d(\gamma _n , U_{C(b)})<1/n$. Let $\gamma _n = c_nz_n$ where $c_n \in U_{C(b)}$, and $z_n \in U_{Z(b)}$. Notice that $z_n \to 1$, $bz_n=z_nb$, and that $b^{-t} c_n b^t \to 1$ as $t \to \infty$.

Choose $t_n >0$ such that $d(b^{-t_n} c_n b^{t_n},1) < 1/n$. Then 
$$ b^{-t_n}\gamma _n b^{t_n} = ( b^{-t_n} c_n b^{t_n})z_n \to 1$$
 By Theorem 1.12 of \cite{Ra 2}, $\{b^{-t}\}_{t>0}$ is not contained in any compact subset of ${\bf G}(\ZZ) \backslash {\bf G}(\RR)$, which contradicts that $b^{-t} \in A_0 \se X_\ZZ$ (Lemma~\ref{l:now}).

\end{proof}

\subsection{$L>0$ and choice of cell in $A_0$.} At this point, we fix $L>0$ to be sufficiently large. We will use this fixed $L$ for our proof of the Theorem~\ref{t:m}.

Let $W_i \se A$ be the kernel of a root $\beta _i \in \Phi _\RR^+$ such that the visual image of $W_i$ in $A^\infty$ is a great sphere that contains $\Sigma _i$.

We let $F$ be the component of $A_0 - \cup _i a_+^L W_i$ that contains $1$.

\begin{lemma} $F$ is compact Euclidean polygon with volume $O(L^{{\text{rk}}_\RR \bf G -1})$.
\end{lemma}

\begin{proof}

The visual cone of $\Sigma$ in $A$ based at $a_+^L$ is a connected component of $A - \cup _i a_+^L W_i$.

The lemma follows if $\Sigma$ and $a_+^\infty$ are contained in distinct components of $A^\infty-A_0^\infty$, and if $a_+^{-\infty} = \lim_{t \to \infty}a_+^{-t} \in \Sigma$. That is indeed the case:  $\alpha (a_+)>1$ for all $\alpha \in \Phi _\QQ (\Delta _\QQ - \alpha _0)^+ $ so  ${\bf P}=  {\bf U}_{\Phi _\QQ (\Delta _\QQ - \alpha _0)^+}{\bf Z_G(T}_{\Delta _\QQ - \alpha _0}) $ fixes $a_+^\infty$. Hence, $a_+^\infty \in \Pi^\infty$. The antipodal map on $A^\infty$ stabilizes $A_0^\infty$, transposes $a_+^{\infty}$ and $a_+^{-\infty}$, and maps $\Pi^\infty$ onto $\Pi_-^\infty \se \Sigma$.

\end{proof}

We denote the face of $F$ given by $a_+^LW_i \cap F$ as $F_i$, so that the topological boundary of $F$ equals $\cup _{i =1}^n F_i$. 

\section{Other cells in $X_\ZZ$ and their homological boundaries}

We denote the real points of the root group ${\bf U}_{(\beta_i)}$ as $U_i$, and  $\langle U_i \rangle _i$ is the group generated by the $U_i$ for $ i \in  \{1,2,...,n\}$.

\begin{lemma}
For each $ i \in  \{1,2,...,n\}$, $U_i \leq U_P$, and thus $\langle U_i \rangle _i \leq U_P$ is abelian.
\end{lemma}

\begin{proof}
Since $\beta _i \in \Phi _\RR^+$, we have $U_i \leq U_Q = U_PU_a$. Either $U_i \leq U_P$ or $U_i \leq U_a \leq {\bf Z}_{\bf G}({\bf T} _{\Delta _\QQ-\alpha _0})$.

Because  ${\bf Z}_{\bf G}({\bf T} _{\Delta _\QQ-\alpha _0})$ is contained in both $\bf P$ and $\bf P^-$, the latter case implies that $U_i$ fixes the antipodal cells $\Pi^\infty$ and $\Pi_-^\infty$. The fixed point set of $U_i$ is a hemisphere in $A^\infty$ with boundary equal to $W_i^\infty$. Thus,  $\Pi^\infty$ and $\Pi_-^\infty$ are contained in $W_i^\infty$, which contradicts that $\Sigma _i = \Sigma \cap W_i^\infty $ does not contain   $\Pi_-^\infty$.

Having ruled out the latter case,  $U_i \leq U_P$ and the lemma follows from Lemma~\ref{l:sa}.

\end{proof}

\subsection{A space for making cycles in $X_\ZZ$.}

\begin{lemma}\label{l:cellhome} $\langle U_i \rangle _i F \se X_\ZZ$.
\end{lemma}

\begin{proof} Because ${\bf R_u(P)}$ is unipotent,  ${\bf R_u(P)} ( \ZZ )$ is a cocompact lattice in $U_P$. We choose a compact fundamental domain $D \se U_P$ for the ${\bf R_u(P)} ( \ZZ )$-action.

There is also a compact set $C \se A_0 = A_\ZZ$  such that ${\bf A}(\ZZ) C =A_\ZZ=A_0$.
As $DC$ is compact, we may assume that ${\bf G} ( \ZZ ) DC  \se X_\ZZ$.

Recall that $\bf A$ is contained in  $\bf P$, so $\bf A$ normalizes ${\bf R_u(P)}$. Hence,
\begin{align*}
 \langle U_i \rangle _i A_0 
 & \se U_P  {\bf A}(\ZZ) C \\
& \se {\bf A}(\ZZ) U_P C \\
& \se {\bf A}(\ZZ) {\bf R_u(P)} ( \ZZ ) D C \\
& \se {\bf G} ( \ZZ )D C \\
& \se X_\ZZ
 \end{align*}
\end{proof}

\subsection{Description of cells used to build our cycle.}

Given $i \in \{1,...,n\}$, let $f_i $ be a point in $F_i$ that minimizes the distance  to $1 \in A$, and let $u_i \in U_i$ be such that $d(u_i f_i, f_i)=1$. Since $F_i \se a_+^LW_i$, any $f \in F_i$ can be expressed as $f=wf_i$ for some $w \in {\emph Ker}(\beta _i)$. It follows that $Ad(w)$ acts trivially on the Lie algebra of $U_i$, that $u_i$ commutes with $w$, and that  $$d(u_i f, f)=d(u_i wf_i, wf_i)=d(wu_i f_i, wf_i)=d(u_i f_i,f_i)=1$$
Setting $\overline{u_i}=\{u_i^t\}_{t=0}^1$, the space $\overline{u_i} F_i$ is a metric direct product of volume $O(L^{{\text{dim}}(F_i)})$.

For $I \subseteq \{1,...,n\}$, let let $F_I =\cap _{i \in I}F_i$ with $F_\emptyset = F$. And let $u_I=\prod _{i \in I} u_i$ and $\overline{u_I}=\prod _{i \in I} \overline{u_i}$ with $\overline{u_\emptyset}={u_\emptyset}=1$.

Similar to the case when $|I|=1$, $\overline{u_I} F_I$ is a metric direct product of volume $O(L^{{\text{dim}}(F_I)})$.

\subsection{Homological boundaries of the cells.}

We endow each interval $\overline{u_i}=[0,u_i]$ with the standard orientation on the closed interval, and we orient  each $\overline{u_I}$ with the product orientation, where the product is taken over ascending order in $\mathbb{N}$. Given $m \in I$, we let $s_I(m)$ be the ordinal of $m$ assigned by the order on $I$ induced by $\mathbb{N}$. Notice that the standard formula for the homological boundary of a cube then becomes $$\partial (\overline{u_I})=\sum _{m \in I} (-1)^{s_I(m)} \big( \overline{u_{I-m}} -  u_m \overline{u_{I-m}} \big)$$

We assign an orientation to $F$, and then assign the orientation to each $F_i$ such that $$\partial (F ) =\sum _{i=1}^nF_i$$

In what follows, if we are given a set $I\se \{1,...,n\}$ with an ordering (which may differ from the standard order on $\mathbb{N}$), and if $m\in \{1,...,n\}$ with $m \notin I$, then the set $I \cup m$ is ordered such that the original order on $I$ is preserved and $m$ is the ``greatest" element of $I \cup m$. For example, $\{1, 7, 5\}\cup 3= \{1,7,5,3\}$.

If $m \in I$, for some ordered set $I\se \{1,...,n\}$, then we endow $I - m$ with the order restricted from $I$.

For an ordered $I$ and $m \in I$, let $r_I(m)=1$ if an even number of transpositions are required to transform the order on $I$ to the order on $(I-m)\cup m$. Let  $r_I(m)=-1$ otherwise.

Given an ordering on a set $I \se \{1,...,n\}$, an orientation on $F_I$, and some $m \in \{1,...,n\}$ with $m \notin I$, we define the orientation of $F_{I\cup m}$ to be such that $F_{I\cup m}$, and not $-F_{I\cup m}$, is the oriented cell that appears as a summand in $\partial( F_I)$. Therefore $$\partial (F_I)= \sum_{m \notin I} F_{I\cup m}$$

In what follows, whenever we write the \emph{exact} symbols $F_I$ or $F_{I'}$ -- but not necessarily the symbol $F_{I\cup m}$ -- the order on $I$ or $I'$ will be the order from $\mathbb{N}$. Thus, the orientation on $F_I$ and $F_{I'}$ can be unambiguously determined from the above paragraph.

It's easy to check that if $I$ is ordered by the standard order on $\mathbb{N}$ and $m \in I$, then $(-1)^{s_I(m)}r_I(m)=(-1)^{|I|}$ and thus $$-(-1)^{s_I(m)}=(-1)^{|I|-1}r_I(m)$$

 Suppose $w_0$ is an outward normal vector for $F_{I\cup m}$ with respect to $F_I$, and $w_1,...w_k$ is a collection of vectors tangent to $F_{I\cup m}$ such that $\{w_0,w_1,...,w_k\}$ defines the orientation for $F_{I}$. Then $\{w_1,...,w_k\}$ defines the orientation for $F_{I\cup m}$. If $\{v_1,...,v_{|I|}\}$ is an ordered basis for the tangent space of $\overline{u_I}$ that induces the standard orientation on $\overline{u_I}$, then $|I|$ transpositions are required to arrange the ordered basis  $$\{w_0,v_{1},...,v_{|I|},w_1,...,w_k\}$$ into the ordered basis   $$\{v_{1},...,v_{|I|},w_0,w_1,...,w_k\}$$ 
 That is, the orientation on $ \overline{u_I}F_{I\cup m}$ defined above is a $(-1)^{|I|}$-multiple of the orientation on $ \overline{u_I}F_{I\cup m}$
 assigned by $\partial ( \overline{u_I}F_{I})$.
 
 It follows from this fact and our above formulas for  $\partial (\overline{u_I})$  and  $\partial (F_I)$ that
  $$\partial (\overline{u_I}F_I)=\sum _{m \in I} (-1)^{s_I(m)} \big( \overline{u_{I-m}} -  u_m \overline{u_{I-m}} \big)F_I +(-1)^{|I|}\sum _{m \notin I} \overline{u_I}F_{I\cup m}$$

\section{A cycle in $X_\ZZ$}

Let $$Y = \sum_{\begin{subarray}{l} K,I \se \{1,...,n\} \\  K\cap I =\emptyset \end{subarray}} (-1)^{|K|} u_K  \overline{u_I} F_I$$

\begin{lemma}
$Y$ is a cycle that is contained in $X _\ZZ$ and has volume $O(L^{{\emph{\text{rk}}}_\RR \bf G -1})$.
\end{lemma}

\begin{proof} Each cell of $Y$ is contained in $X _\ZZ$ by Lemma~\ref{l:cellhome} and has volume $O(L^k)$ for $k \leq {\text{rk}}_\RR {\bf G} -1$, so we have to check that $\partial Y =0$.

From our formula for  $\partial (\overline{u_I}F_I)$ we have that
\begin{align*} 
\partial Y & = \sum_{\begin{subarray}{l} K,I \se \{1,...,n\} \\  K\cap I =\emptyset \end{subarray}} (-1)^{|K|} u_K  \Big[\sum _{m \in I} (-1)^{s_I(m)} \big( \overline{u_{I-m}} -  u_m \overline{u_{I-m}} \big)F_I   \\ 
& \hspace{.5in} +(-1)^{|I|}\sum _{m \notin I} \overline{u_I}F_{I\cup m} \Big] \\
& = \sum_{\begin{subarray}{l} K,I \se \{1,...,n\} \\  K\cap I =\emptyset \end{subarray}}  \sum _{m \in I} (-1)^{s_I(m)} (-1)^{|K|} u_K \big( \overline{u_{I-m}} -  u_m \overline{u_{I-m}} \big)F_I   \\ 
& \hspace{.5in} + \sum_{\begin{subarray}{l} K,I \se \{1,...,n\} \\  K\cap I =\emptyset \end{subarray}} (-1)^{|I|}\sum _{m \notin I} (-1)^{|K|} u_K \overline{u_I}F_{I\cup m} 
\end{align*}

For $K,I \se \{1,...,n\} $ with $  K\cap I =\emptyset$ we have

\begin{align*} \sum_{m \notin I}(-1)^{|K|}   &  u_K\overline{u_I} F_{I\cup m} \\ 
&  = \sum_{m \notin I\cup K} (-1)^{|K|}   u_K\overline{u_I} F_{I\cup m} \\
&  \hspace{.5in}  + \sum_{m \in K }(-1)^{|K|}   u_K\overline{u_I} F_{I\cup m} \\
&  = \sum_{m \notin I\cup K}(-1)^{|K|}   u_K\overline{u_{(I\cup m)-m}} F_{I\cup m}
\\ &  \hspace{.5in} + \sum_{m \in K }(-1)^{|K|}    u_{K-m} u_m\overline{u_{(I\cup m)-m}} F_{I\cup m} \end{align*}

There is a natural bijection between triples $(I,K,m)$ where  $K\cap I =\emptyset$ and $m \notin I \cup K$, and triples $(I',K',m)$ where $K'\cap I' =\emptyset$ and $m \in I'$. To realize the bijection, let $K'=K=K-m$ and $I'=I\cup m$.

There is also a bijection between triples $(I,K,m)$ where  $K\cap I =\emptyset$ and $m \in K$, and triples 
$(I',K',m)$ where $K'\cap I' =\emptyset$ and $m \in I'$. This bijection is also realized by setting $K'=K-m$ and $I'=I\cup m$.

Therefore, if we let $K'=K-m$ and $I'=I\cup m$ then the above equation gives 

\begin{align*} \sum_{\begin{subarray}{l} K,I \se \{1,...,n\} \\  K\cap I =\emptyset \end{subarray}}  & (-1)^{|I|} \sum_{m \notin I}(-1)^{|K|}    u_K\overline{u_I} F_{I\cup m}  \\
& =\sum_{\begin{subarray}{l} K',I' \se \{1,...,n\} \\  K'\cap I' =\emptyset \end{subarray}}  (-1)^{|I'|-1}  \Big[ \sum _{m \in I' } (-1)^{|K'|}  r_{I'}(m) u_{K'}\overline{u_{I'-m}} F_{I'}  \\
& \hspace{.5in} +\sum_{m \in I' } (-1)^{|K'\cup m|}  r_{I'}(m)  u_{K'} u_m \overline{u_{I'-m}} F_{I'}\Big] \\
& =\sum_{\begin{subarray}{l} K',I' \se \{1,...,n\} \\  K'\cap I' =\emptyset \end{subarray}}  (-1)^{|I'|-1} \Big[ \sum _{m \in I' } (-1)^{|K'|}  r_{I'}(m) u_{K'}\overline{u_{I'-m}} F_{I'}  \newline \\
& \hspace{.5in} - \sum _{m \in I' } (-1)^{|K'|}  r_{I'}(m) u_{K' }u_m\overline{u_{I'-m}} F_{I'}\Big] \\
& =\sum_{\begin{subarray}{l} K,I \se \{1,...,n\} \\  K\cap I =\emptyset \end{subarray}}  (-1)^{|I|-1}  \Big[ \sum _{m \in I} (-1)^{|K|}  r_{I}(m) u_{K}\overline{u_{I-m}} F_{I}  \newline \\
& \hspace{.5in} - \sum _{m \in I } (-1)^{|K|}  r_{I}(m) u_{K }u_m\overline{u_{I-m}} F_{I}\Big] \\
& =\sum_{\begin{subarray}{l} K,I \se \{1,...,n\} \\  K\cap I =\emptyset \end{subarray}}  (-1)^{|I|-1}   \sum _{m \in I} (-1)^{|K|}  r_{I}(m)  u_{K} \big(\overline{u_{I-m}} -  u_{ m} \overline{u_{I-m}} \big) F_{I} \\
& =\sum_{\begin{subarray}{l} K,I \se \{1,...,n\} \\  K\cap I =\emptyset \end{subarray}}     \sum _{m \in I} (-1)^{|I|-1} r_{I}(m)  (-1)^{|K|}   u_{K} \big(\overline{u_{I-m}} -  u_{ m} \overline{u_{I-m}} \big) F_{I} \\
& =- \sum_{\begin{subarray}{l} K,I \se \{1,...,n\} \\  K\cap I =\emptyset \end{subarray}}    \sum _{m \in I} (-1)^{ s_{I}(m)}  (-1)^{|K|}   u_{K} \big(\overline{u_{I-m}} -  u_{ m} \overline{u_{I-m}} \big) F_{I}
  \end{align*} 

Substituting the preceding equation into our equation for $\partial Y$ proves
$$\partial Y =0 $$

\end{proof}

\section{Fillings of Y}

There exists polynomially efficient fillings for $Y$ in the symmetric space $X$.

\begin{lemma}\label{l:poly} There exists a chain $Z$ with volume $O(L^{{\emph{\text{rk}}}_\RR \bf G })$ and  $\partial Z =Y$.
\end{lemma}
\begin{proof}

As $Y \se \overline{u_I} F$, it follows from Lemma~\ref{l:ad} that there is some $T=O(L)$ such that $a_+^T Y$ is contained in an $\varepsilon$-neighborhood of $a_+^TF$, which is isometric to $F$.  Thus, there is a filling, $Z_0$, of $a^T_+Y$ of volume $O(L^{{\text{rk}}_\RR \bf G -1})$.

Let $Z=Z_0 \cup _{t\in \{1,T\}}a_+^TY$.

\end{proof}

\subsection{Fillings of $Y$ in $X_\ZZ$.} In contrast to Lemma~\ref{l:poly}, the fillings of $Y$ that are contained in $X_\ZZ$ have volumes bounded below by an exponential in $L$. A fact that we will prove after a couple of helpful lemmas.

For $f \in F$, define $d_i(f)$ to be the distance in the flat $A$ between $f$ and $a_+^LW_i$. 

\begin{lemma}
There are $s_i>1$ and $ s_0 >0$ such that the cube  $\overline{u_I} f$ with the path metric is isometric to $\prod _ {i \in I} [0,e^{s_id_i(f)+s_0}]$.
\end{lemma}

\begin{proof}
It suffices to prove that $\overline{u_i} f$ is isometric to $ [0,e^{s_id_i(f)+s_0}]$.

Choose $b_i \in A$ such that $d(b_i, 1)=d(f,a_+^LW_i)=d_i(f)$ and such that there exists some $w_i \in W_i$ with $f=b_ia_+^Lw_i$. Notice that  $W_i$ separates $b_i$ from $a_+^L$ in $A$. Since $U_i \leq U_P$, Lemma~\ref{l:ad} shows that $\beta_i(a_+^L)>1$. It follows that $\beta _i (b_i)<1$.

With $d_\Omega$ as the path metric of a subspace $\Omega \se X$,
$$d_{U_if}(u_if,f)=d_{U_if}(u_ib_ia_+^Lw_i,b_ia_+^Lw_i) $$ 

As $W_i$ is the kernel of $\beta _i$, $w_i$ commutes with $u_i$ implying
\begin{align*} d_{U_if}(u_if,f) & =d_{w_i^{-1}U_if}(u_ib_ia_+^L,b_ia_+^L) \\
&= d_{U_i}(a_+^{-L} b_i^{-1}u_ib_i a_+^L,1) \end{align*}

On the Lie algebra of $U_i$, $Ad(a_+^{-L} b_i^{-1})$ scales by $\beta_i(a_+^{-L})\beta _i( b_i)^{-1}$.

\end{proof}

In the above lemma we may let $f=1$ and let $I$ be the singleton $i$. It can easily be seen that $d_i(1)=O(L)$ which leaves us

\begin{lemma} There is some $C>0$ such that  $d_{U_i}(u_i,1) \geq  e^{C L+s_0}$ for any $i$.

\end{lemma}

We conclude our proof of Theorem~\ref{t:m} with the following

\begin{lemma}\label{l:exp}
Suppose there is a chain $B \se X_\ZZ$ such that $\partial B=Y$. Then the volume of $B$ is bounded below by $e^{C_0L}$ for some $C_0>0$.
\end{lemma}

\begin{proof}
Suppose $B$ has volume $\lambda$. By Lemma~\ref{l:push}, $\pi(B) \se  U_QA_0$ has volume $O(\lambda)$.

Recall that $Y \se U_QA_0$, so $\partial \pi (B) =Y$.

After perturbing $\pi (B)$, we may assume that $\pi (B)$ is transverse to $U_Q$, and that the $1$-manifold $\pi (B) \cap U_Q$ has length proportional to the volume of $\pi(B)$. Since $$\partial ( \pi (B) \cap U_Q)=\partial \pi (B) \cap U_Q = Y \cap U_Q = \{u_I\}_{I\se \{1,...,n\}}$$  there is an $I\se \{1,...,n\}$ and a path $\rho : [0,1] \rightarrow \pi (B) \cap U_Q$ such that $\rho (0)=1$ and $\rho(1)=u_I $ with $\text{length}(\rho ) = O(\lambda)$.

Choose $i \in I$. $U_Q$ is nilpotent, so the distortion of the projection $q: U_Q \rightarrow U_i$ is at most polynomial. Therefore, $q \circ \rho$ is a path in $U_i$ between $1$ and $u_i$ with $\text{length}(q \circ \rho)=O(\lambda^k)$ for some $k \in \mathbb{N}$.

The preceding lemma showed  $e^{C L+s_0} \leq  \text{length}(q \circ \rho)$. Therefore, $\lambda \geq \kappa e^{\frac{C}{k} L}$ for some $\kappa >0$.

\end{proof}

Combining Lemmas~\ref{l:poly} and~\ref{l:exp} yields Theorem~\ref{t:m}.

\end{document}